\documentclass[12pt]{article}

\usepackage[cp1251]{inputenc}
\usepackage[english]{babel}

\usepackage{amsfonts,epsfig,epstopdf}
\usepackage{amssymb,amsmath,enumerate,float}

\usepackage{graphicx}
\usepackage{subfigure}
\usepackage{color}
\definecolor{light-gray}{gray}{0.6}
\definecolor{dark-gray}{gray}{0.35}
\usepackage[bookmarksopen=true,bookmarks=true,bookmarksnumbered=true,colorlinks=true,linkcolor=dark-gray,citecolor=light-gray]{hyperref}
\usepackage{mathrsfs}
\usepackage{authblk}

\makeatother

\addto\captionsenglish{}

\newcommand{\be}{\begin{equation}}
\newcommand{\ee}{\end{equation}}
\newcommand{\beq}{\begin{equation}}
\newcommand{\eeq}{\end{equation}}

\newcommand{\rin}{{\rm int}}
\newcommand{\rex}{{\rm ext}}
\newcommand{\K}{\mathcal{K}}
\newcommand{\A}{\mathbf{A}}
\newcommand{\B}{\mathbf{B}}
\newcommand{\C}{\mathbf{C}}
\newcommand{\matPI}{\mathbf{\Pi}}
\newcommand{\Zero}{\mathbf{0}}
\newcommand{\Id}{\mathbf{I}}

\newcommand{\hfrac}[2]{#1 / #2} 

\newcommand{\fref}[1]{Fig.~\ref{#1}}

\newcommand{\sref}[1]{Section~\ref{#1}}

\newcommand{\cref}[1]{Chapter~\ref{#1}}

\newcommand{\aref}[1]{Appendix~\ref{#1}}
\providecommand{\keywords}[1]{\textbf{\textit{Keywords:~}} #1}

\title{A New Numerical Method for Solving the Acoustic Radiation Problem \footnote{The paper is based on the talk given at the Fourth International Symposium and the Sixth All-Russia Conference `Computational Experiment in Aeroacoustics' held on September 19-24, 2016, in Svetlogorsk, Kaliningrad region, http://ceaa-w.ima-mod.ru. The present document includes some more numerical simulations with respect to the journal version.}}

\author[1]{J. Poblet-Puig\thanks{correspondence: UPC, Campus Nord B1, Jordi Girona 1, E-08034 Barcelona, Spain, e-mail: jordi.poblet@upc.edu}}
\author[2]{A. V. Shanin\thanks{e-mail: a.v.shanin@gmail.com}}
\affil[1]{Laboratori de C\`{a}lcul Num\`{e}ric, E.T.S. d'Enginyers de Camins, Canals i Ports de Barcelona, Universitat Polit\`{e}cnica de Catalunya}
\affil[2]{Department of Physics, Acoustics Division, Moscow State University}
\begin{document}

\maketitle

\begin{abstract}



A numerical method of solving the problem of acoustic wave radiation in the presence of a rigid
scatterer is described. It combines the finite element method and the boundary algebraic equations. In the
proposed method, the exterior domain around the scatterer is discretized, so that there appear an infinite
domain with regular discretization and a relatively small layer with irregular mesh. For the infinite regular
mesh, the boundary algebraic equation method is used with spurious resonance suppression according to
Burton and Miller. In the thin layer with irregular mesh, the finite element method is used. The proposed
method is characterized by simple implementation, fair accuracy, and absence of spurious resonances.

\end{abstract}

\section*{List of symbols and acronyms}

\begin{tabular}{ll}
$\alpha,\beta$ & coefficients of the numerical technique considered in each domain\\
$\A,\B,\C$ & matrices of the method \\
$\matPI$ & projector matrix \\
$BAE$ & Boundary Algebraic Equations \\
$BEM$ & Boundary Element Method \\
$CFIE$ & Combined-Field Integral Equations \\
$DtN$ & Dirichlet-to-Neumann\\
$\delta_{j,m}$ & Dirac delta \\
$FEM$ & Finite Element Method \\
$f$ & force term, sources of the field \\
$G_{j,m}$ & discrete Green's function \\
$\Gamma_\rin$ & scaterer surface (or curve) \\
$\Gamma_\rex$ & boundary between domains $\Omega_\rin$ (solved with FEM) \\
              & and $\Omega_\rex$ (solved with BAE) \\
$\gamma_\rex$ & set of nodes on $\Gamma_\rex$\\
$\gamma_{\rm o}$ & set of nodes surrounding $\gamma_\rex$ and $\gamma_\rex$ itself\\
$h$ & grid or finite element size \\
$h^\rin, h^\rex$ & fluxes across $\Gamma_\rex$\\
$\K$ & wavenumber \\
$R$ & radius of the circular scatterer \\
$u$ & main variable (scattered field) \\
$\Omega_\rin$ & domain inside $\Gamma_\rex$ and around the scatterer \\
$\omega_\rin$ & set of nodes in $\Omega_\rin$ \\
$\omega'_\rin$  & set of elements in $\Omega_\rin$ \\
$\Omega_\rex$ & infinite domain outside $\Gamma_\rex$ \\
$\omega_\rin$ & set of nodes in $\Omega_\rex$ \\
$\omega'_\rin$  & set of elements in $\Omega_\rex$ \\
$\Omega$ & entire space covered with uniform (periodic) mesh  \\
$\omega$ & set of nodes in $\Omega$ \\
$\omega'$  & set of elements in $\Omega$ \\
\end{tabular}

\section{Introduction}

The problem of external acoustic scattering has recently been solved \cite{poblet-PVS:2015} by means of the boundary algebraic equations method (BAE \cite{Martinsson-Rodin:2009,Gillman-Martinsson:2010,Tsukerman:2011,Bhat-Osting:2009}) and considering a combined-field integral formulation (CFIE, \cite{Burton-Miller:1971,kirkup:1998}). The resulting method is, essentially, a discrete analogue of the boundary element method (BEM) that inherits the good properties of BAE and the advantages of CFIE, avoiding most of the BEM drawbacks. On the one hand, no quadratures of oversingular integrals have to be computed due to the discrete nature of BAE. On the other hand, the resulting integral equations are free of spurious resonances due to the CFIE formulation \cite{poblet-PVS:2014}.

However, the main drawback of the CFIE--BAE method is the reduction of accuracy when smooth scatterers with curved surfaces such as spheres or cylinders are considered. This is because the method is based on a regular discretisation of the space (grid) and the obstacles must be approximated by means of the closest brick-description. The proper description of arbitrary shaped scatterers is a common problem in the family of methods based on regular grids, see for example \cite{medvinsky-MTT:2013}.

Our goal here is to present a complementary formulation where the CFIE--BAE method is coupled with some more versatile numerical technique in order to deal with arbitrary shaped scatterers. This will typically be a thin layer of finite elements (FEM) between the obstacle surface and a close grid-shaped boundary that surrounds the obstacle. The FEM domain has on the one side the boundary conditions corresponding to the scatterer and on the other side the coupling with the CFIE--BAE. This acts as a method for domain truncation and exactly imposes the radiation boundary conditions.

The coupling of numerical methods in order to maximize the benefits and reduce the disadvantages of each one has been often used. See for example \cite{Zienkiewicz-ZKB:1977} where the FEM was complemented with a boundary integral method to deal with radiation conditions, \cite{JohnsonNedelec:1980} where the stability conditions of FEM--BEM couplings were studied.

Some more recent works on the FEM--BEM coupling applied to the scattering of waves can be found, see for example \cite{hsiao:1991, Chiang:2000, Gatica:2009}. However, to the best of the authors knowledge, the coupling BAE--FEM has not been considered.

The method presented here can also be understood as an alternative to impose the radiation boundary condition and truncate the computation of domains. It has the added value that the obtained solutions are `exact' in the sense that no numerical artefact is required. In some popular alternatives such as the perfectly matched layers (PML \cite{Berenger:1994}) the reflected waves are attenuated by means of a virtual damping medium placed in the surrounding of the problem domain. It certainly diminishes the reflected waves but it is well know that their parameters (i.e. complex wave number of the medium) must be calibrated properly. Moreover, evanescent waves can remain undamped (see \cite{zampolli:2008,basu:2003}) and the quality of solution can be diminished in some zones close to the layer such as the corners. The shape of the PML, the thickness of the layer and the distance from the scatterer are important aspects also for the quality of the solution and in order to derive the PML equations. On the contrary, the approach presented here is more flexible in the sense that it is independent of the shape and the outer boundary can be placed very close to the scatterer without affecting the quality of the solution. This will be illustrated later in \sref{sec:NumericalExamples}.

In the remainder of the document, the formulation of the problem  is presented in \sref{sec:formulation} and the method is detailed in \sref{sec:method}. Its properties are shown with the numerical examples in \sref{sec:NumericalExamples} before the conclusions. The parts of the development that are not essential have been grouped in the appendices: some details of the derivation of BAE equations in \aref{sec:AppA} and a proof of solvability in \aref{sec:AppB}.


\section{Formulation of continuous and discrete problems}
\label{sec:formulation}

We consider a 2D of 3D external acoustic stationary problem. The scatterer is approximated by a surface (or a curve)
$\Gamma_\rin$.
The inhomogeneous Helmholtz equation
\begin{equation}
\Delta u + \K^2 u = f
\label{eq0201}
\end{equation}
is assumed to be fulfilled in the medium. Variable $u$ may correspond to acoustical pressure or acoustical potential.
We assume that the boundary is acoustically hard (Neumann).

\begin{figure}[ht]
\centerline{\epsfig{file=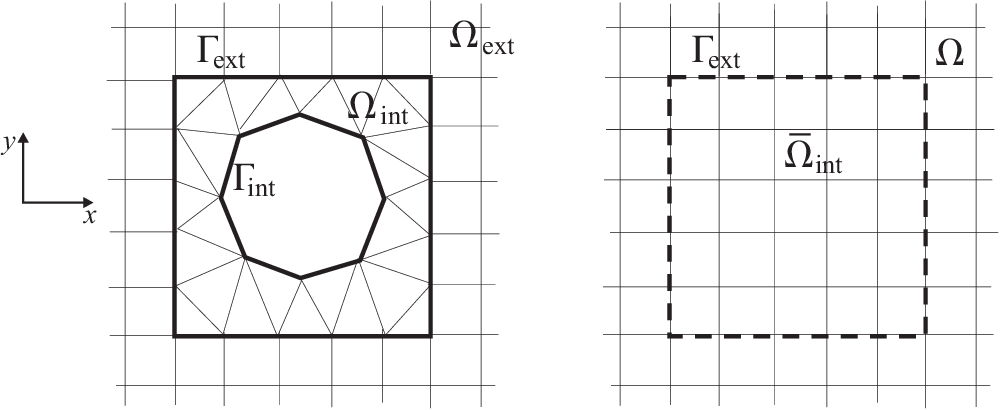,width = 12cm}}
\caption{Two domains and meshes (left), the uniform mesh (right)}
\label{fig01}
\end{figure}

Function $f$ represents the sources of the field, i.\ e.\  a radiation problem is studied. If the sources are
put on the surface $\Gamma_\rin$ then one can study radiation of wave by a vibrating boundary. Typically it is necessary to find directivity of the field as the result.

We assume that the exponential factor of an outgoing wave has form of $\exp \{ i \K r\}$ for big $r$, where $r$ is the distance from the origin.
We assume that $\K$ has a vanishing positive imaginary part. Thus, an outgoing wave should decay
exponentially at infinity. The same property (it is the {\em radiation condition}) should be obeyed by
any numerical approximation of~$u$.


Split the domain external with respect to $\Gamma_\rin$ into two subdomains $\Omega_\rin$ and $\Omega_\rex$
(one inside another, see \fref{fig01}). The boundary between these domains, $\Gamma_\rex$, should have a simple shape. For example, the interior of $\Gamma_\rex$ should be  a union of equal cubes/squares. This property will enable us to apply the BAE method to $\Gamma_\rex$.

The boundary $\Gamma_\rex$ does not correspond to any physical interface, but it divides the space into two parts, which be treated numerically in a different manner. The wave process in domain $\Omega_\rex$ will be modelled
by the BAE method, thus giving a boundary condition (an approximation of a DtN operator) on $\Gamma_\rex$.
This boundary condition should establish the absence of waves coming from infinity on $\Gamma_\rex$.
The internal domain $\Omega_\rin$ will be treated by a usual finite element method.

We assume that all sources lie inside $\Gamma_\rex$, i.\ e.\ belong to $\Omega_\rin$.

Consider the entire space $\Omega$ covered with uniform (periodic) mesh (see \fref{fig01}, right). The mesh in our understanding consists of nodes and finite elements (polygons or polyhedra). The set of all nodes belonging to the uniform mesh will be denoted by~$\omega$, and the set of all finite elements by $\omega'$.
Let $\bar \Omega_\rin$ be a domain composed of some
finite elements of the uniform mesh $\bar \omega_\rin'$. Denote the set of nodes adjacent to these selected elements by $\bar \omega_\rin$. Denote the set of nodes adjacent to the finite elements
$\omega_\rex' = \omega' \setminus \bar \omega'_\rin$ of the uniform mesh by $\omega_\rex$. The boundary nodes form the
set $\gamma_\rex = \bar \omega_\rin \cap \omega_\rex$. Obviously, these nodes belong to $\Gamma_\rex$.

Now consider a non-uniform mesh defined in domain $\Omega_\rin$ (\fref{fig01}, left).
Denote the set of nodes of this mesh by
$\omega_\rin$ and the elements of this mesh by $\omega'_\rin$.
The nodes of $\omega_\rin$ lying on the boundary $\Gamma_\rex$
should coincide with $\gamma_\rex$, i.\ e.\ the uniform mesh on $\Omega_\rex$ and the arbitrary mesh
on $\Omega_\rin$ should form together a valid mesh on $\Omega_\rin \cup \Omega_\rex$.
Also, $\omega_\rin \cap \omega_\rex = \gamma_\rex$.

Let equation
\begin{equation}
\Delta u + \K^2 u = g,
\label{eq0201a}
\end{equation}
valid in the entire space $\Omega$,
be approximated on a uniform mesh $\omega, \omega'$ using the finite element
method. Let the nodal values of $u$ and $g$ be denoted by $u_j$, $g_j$. Write the approximation
in the form
\begin{equation}
\sum_{k \in \omega} \beta_{j,k} u_k = g_j ,
\qquad
j \in \omega .
\label{eq0202}
\end{equation}
Assume that the coefficients $\beta_{j,k}$ possess the following properties:
\begin{itemize}
\item
$\beta_{j,k} \ne 0$ only for nodes $j$ and $k$ adjacent to the same finite element;
\item
the matrix is symmetrical $\beta_{j,k} = \beta_{k,j}$ ;
\item
since the mesh is periodical, the coefficients do not change when the pair of nodes is translated
along the mesh.
\end{itemize}

Now consider the approximation of equation (\ref{eq0201}) in the domain $\Omega_\rin \cup \Omega_\rex$. Let this
approximation be written in the form
\begin{equation}
\sum_{k \in (\omega_\rin \cup \omega_\rex) }\alpha_{j,k} u_k
= f_j, \qquad j \in \omega_\rin \cup  \omega_\rex ,
\label{eq0203}
\end{equation}
Let the coefficients $\alpha_{j,k}$ have the following properties:
\begin{itemize}
\item
$\alpha_{j,k} \ne 0$ only for nodes $j$ and $k$ both belonging to the same finite element;
\item
the matrix is symmetrical $\alpha_{j,k} = \alpha_{k,j}$ ;
\item
$\alpha_{j,k} = \beta_{j,k}$ if $j, k \in \omega_\rex$ and  at least one of the nodes $j$, $k$ belongs to
$\omega_\rex \setminus \gamma_\rex$.
\end{itemize}

The last point means that the discretisation (\ref{eq0203})
is uniform in $\Omega_\rex$.

Since the Neumann boundary condition is imposed on $\Gamma_\rin$, equation (\ref{eq0203}) naturally incorporates the
boundary condition. The method, though, can be easily modified to the case of arbitrary boundary conditions.

Our aim is to present a method for solving (\ref{eq0203}). Equation (\ref{eq0202}) is auxiliary for the method.


\section{FEM--BAE method}
\label{sec:method}

Split equation (\ref{eq0203}) into two equations:
\begin{equation}
\sum_{k \in \omega_\rin} \alpha_{j,k}^\rin u_{k}^\rin = f_j + h_j^\rin,
\qquad
j \in \omega_\rin
\label{eq0301}
\end{equation}
\begin{equation}
\sum_{k \in \omega_\rex} \alpha_{j,k}^\rex u_{k}^\rex =  h_j^\rex,
\qquad
j \in \omega_\rex.
\label{eq0302}
\end{equation}
The matrices $\alpha^\rin_{j,k}$, $\alpha^\rex_{j,k}$ and the flows $h^\rex_j$,
$h^\rin_j$ should posses the following properties:
\begin{itemize}
\item
$\alpha_{j,k}^\rin = \alpha_{j,k}$ if $j, k \in \omega_\rin$, and at least one of the nodes
$j, k$ belongs to $\omega_\rin \setminus \gamma_\rex$;
\item
$\alpha_{j,k}^\rex = \alpha_{j,k} = \beta_{j,k}$ if $j, k \in \omega_\rex$, and at least one of the nodes
$j, k$ belongs to $\omega_\rex \setminus \gamma_\rex$;
\item
$\alpha_{j,k}^\rex + \alpha_{j,k}^\rin = \alpha_{j,k}$ if $j, k \in \gamma_\rex$ ;
\item
matrices are symmetrical: $\alpha^\rex_{j,k} = \alpha^\rex_{k,j}$,
$\alpha^\rin_{j,k} = \alpha^\rin_{k,j}$;
\item
$h_j^\rex \ne 0$ or $h_j^\rin \ne 0$ only if $j \in \gamma_\rex$
\item
$h^\rex_j = - h^\rin_j$ if $j \in \gamma_\rex$.
\end{itemize}

Matrices $\alpha_{j,k}^\rex$ and $\alpha_{j,k}^\rin$ possessing the listed properties can be obtained by
assembling the standard FEM matrices performing summation only over the elements belonging to $\omega_\rex'$
or over $\omega_\rin'$, respectively. The flows $h^\rex_j$,
$h^\rin_j$ remain unknown at this stage.

Let also be $u^\rex_j = u^\rin_j$ for $j \in \gamma_\rex$.

By summing (\ref{eq0301}) and (\ref{eq0302})
it is easy to check that the function
\begin{equation}
u_j = \left\{ \begin{array}{ll}
u_j^\rex & j \in \omega_\rex \\
u_j^\rin & j \in \omega_\rin
\end{array}\right.
\label{eq0303}
\end{equation}
is a solution of (\ref{eq0203}). Our plan is to substitute (\ref{eq0302}) by a relation of the form
\begin{equation}
h^\rex_j = \sum_{k \in \gamma_\rex} B_{j,k} u^\rex_k ,
\qquad
j \in \gamma_\rex
\label{eq0304}
\end{equation}
for some matrix $B$, and then represent (\ref{eq0301}) in the form
\begin{equation}
\left( \alpha^\rin_{j,k} + \sum_{m,n \in \gamma_\rex}\Pi^T_{j,m}B_{m,n}\Pi_{n,k} \right) u_k = f_j,
\label{eq0305}
\end{equation}
where $\Pi_{m,n}$, $m \in \gamma_\rex$,  $n \in \omega_\rin$  is a projector matrix
\begin{equation}
\Pi_{m,n} = \left\{   \begin{array}{ll}
1, &  m=n,\quad n\in \gamma_\rex \\
0, &  \mbox{otherwise}
\end{array} \right.
\label{eq0306}
\end{equation}
and $\Pi^T_{m,n} = \Pi_{n,m}$. Then (\ref{eq0305}) can be solved as a linear system.

Expression (\ref{eq0304}) can be obtained from the BAE--CFIE method \cite{poblet-PVS:2015}.
Here we follow the consideration of \cite{poblet-PVS:2015}.
Let $G_{m,n}$ be an approximation of the Green's function of equation (\ref{eq0201a}), i.\ e.\ let
$G_{m,n}$
obey
equation
\begin{equation}
\sum_{k \in \omega}\beta_{j,k} G_{k,m} = \delta_{j,m},
\qquad
j,m \in \omega,
\label{eq0307}
\end{equation}
and the radiation condition. Here $\delta_{j, m}$ is the Kronecker's delta.
Since (\ref{eq0201a}) is an equation on a uniform (periodic) mesh covering the whole space, function $G$ can be computed analytically by the Fourier transformation method. Matrix $G_{m,n}$ is symmetrical:
$G_{m,n} = G_{n,m}$ (see \cite{poblet-PVS:2015}).
Introduce a notation
\begin{equation}
b_{j,m} =
\sum_{n \in \omega_\rex}
\alpha^\rex_{j,n} G_{n,m} - \delta_{j,m}  ,
\qquad
j,m \in \omega_\rex.
\label{eq0308a}
\end{equation}
where $b_{j,m} \ne 0$ only if $j \in \gamma_\rex$ (note that for $j \in (\omega_\rex \setminus \gamma_\rex)$
$\alpha^\rex_{j,n} = \beta_{j,n}$, and (\ref{eq0307}) can be applied).

According to \cite{poblet-PVS:2015}, the BAE--CFIE equation connecting $h^\rex_j$ and $u^\rex_j$, $j \in \gamma_\rex$ is as
follows:
\begin{equation}
\sum_{j \in \gamma_\rex} u^\rex_j A_{j,m} =
\sum_{j \in \gamma_\rex} h^\rex_j
C_{j,m},
\qquad
j,m \in \gamma_\rex,
\label{eq0308b}
\end{equation}
\begin{equation}
A_{j,m} = \delta_{j,m} + b_{j,m} + \nu \sum_{n \in \omega_\rex} b_{j,n} \alpha_{n,m}^\rex ,
\label{eq0309}
\end{equation}
\begin{equation}
C_{j,m} =
- \nu \delta_{j,m} + G_{j,m} + \nu \sum_{n \in \omega_\rex}G_{j,n} \alpha_{n,m}^\rex.
\label{eq0310}
\end{equation}
$\nu$ is an arbitrary complex number with a non-zero imaginary part.

It follows from (\ref{eq0308b}) that matrix $B$ from (\ref{eq0304}) can be written as
\begin{equation}
\B = (\A \C^{-1})^T.
\label{eq0311}
\end{equation}

A known problem associated with the boundary integral equation is linked with formula (\ref{eq0311}) or a similar one. Although $\B$ should exist for all temporal frequencies, if no special measures are undertaken
matrices $\A$ and $\C$ may be singular. This feature is named spurious resonances. For example,
if $\nu =0$ (\ref{eq0308a}) corresponds to Kirchhoff formulation of boundary integral equations.
The Kirchhoff boundary integral equations are known to be prone to spurious resonances \cite{Schenck:41,Benthien-Schenck:1997,chien:1990}.
The CFIE approach is necessary to suppress the spurious resonances. The case ${\rm Im}[\nu] \ne 0$
corresponds to a CFIE formulation.

A sketch of derivation of (\ref{eq0308b}) and a proof of invertibility of $\C$ under some general condition
can be found in the Appendix.

Introduce the set of nodes $\gamma_{\rm o}$ belonging to $\omega_\rex$ and neighbouring $\gamma_\rex$
(i.~e.\ they are the nodes adjacent to the finite elements adjacent to nodes from $\gamma_\rex$).
The set $\gamma_{\rm o}$ is finite. By construction, $\gamma_\rex \subset \gamma_{\rm o}$. The summation
in (\ref{eq0309}) and (\ref{eq0310}) can be held along $\gamma_{\rm o}$ instead of $\omega_\rex$.

Let us summarize the procedure of solving (\ref{eq0203}).
\begin{itemize}
\item
The Green's function $G_{m,n}$ and values $b_{m,n}$ should be tabulated for
$m\in \gamma_\rex$, $n \in \gamma_{\rm o}$.
\item
Matrices $\A$, $\C$ should be calculated from (\ref{eq0309}), (\ref{eq0310}) for $j,m \in \gamma_\rex$.
\item
Matrix $\B$ should be found from (\ref{eq0311}).
\item
Equation (\ref{eq0305}) should be solved.
\end{itemize}

As the result of this procedure, one obtains the nodal values of field $u_j^\rin$. Thus,
the near field becomes known. To get the far field, one needs to perform an additional step of post-processing.
Namely, for any $m \in \omega_\rin$
\begin{equation}
u_m^\rex = \sum_{j \in \gamma_\rex} (h_j^\rex G_{j,m} - u_j^\rin b_{j,m}).
\label{eq0312a}
\end{equation}
Substituting (\ref{eq0304}), obtain
\begin{equation}
u_m^\rex = \sum_{j \in \gamma_\rex} u_j^\rin
\left(
\sum_{k\in \gamma_\rex} B_{k,j} G_{k,m} -   b_{j,m}
\right).
\label{eq0312b}
\end{equation}
If node $m$ is located far enough, asymptotic expressions for $G_{j,m}$ and $b_{m,j}$ can be found.
Formula (\ref{eq0312b}) provides the solution in the far field (a directivity can be taken from it).

It can be convenient to solve the whole problem at the same time and avoid the explicit inversion of matrix $\C$.
One should consider a linear system of equations where the unknowns are $\mathbf{u}^\rex$ and $\mathbf{u}^\rin$ that contain the nodal values in $\gamma_\rex$ and $\omega_\rin$ respectively, and $\mathbf{h}^\rex$ that contain the fluxes $h^\rex$ defined in (\ref{eq0302}). The coupled linear system of equations is
\begin{equation}\label{eq:CoupledSystem}
\begin{bmatrix}
    \A   & \Zero & -\C \\[2.5ex]
    \Zero  & \A^\rin & \matPI^T\\[2.5ex]
    \Id  & -\matPI & \Zero
\end{bmatrix}
\begin{bmatrix}
 \mathbf{u}^\rex \\[2.5ex]
 \mathbf{u}^\rin \\[2.5ex]
 \mathbf{h}^\rex
\end{bmatrix} =
\begin{bmatrix}
 \Zero \\[2.5ex]
 \mathbf{f} \\[2.5ex]
 \Zero
\end{bmatrix}
\end{equation}
where $\A$ and $\C$ are the matrices defined in (\ref{eq0309}) and  (\ref{eq0310}), $\A^\rin$ is the matrix obtained from (\ref{eq0301}) which is typically the usual FEM matrix, $\Zero$ is a null matrix, $\Id$ the identity and $\matPI$ the projector matrix defined in (\ref{eq0306}) (rows for the nodes in $\gamma_\rex$ and columns for the nodes in $\omega_\rin$). The force vector includes $\mathbf{f}$ from (\ref{eq0301}).

In the linear system (\ref{eq:CoupledSystem}) the first block of equations represent (\ref{eq0308b}), the second block of equations accounts for (\ref{eq0301}) and the continuity of fluxes $h^\rex_j = - h^\rin_j$ if $j \in \gamma_\rex$. And finally the third block imposes continuity of variable $u$: $u^\rex_j = u^\rin_j$ for $j \in \gamma_\rex$.



\section{Numerical results}
\label{sec:NumericalExamples}

\begin{figure}[ht]
\subfigure[]{\includegraphics[width=0.4\textwidth]{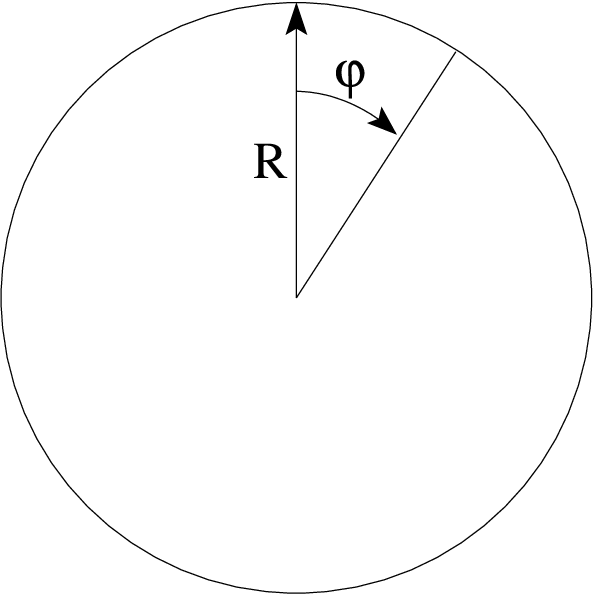}} \subfigure[]{\includegraphics[width=0.45\textwidth]{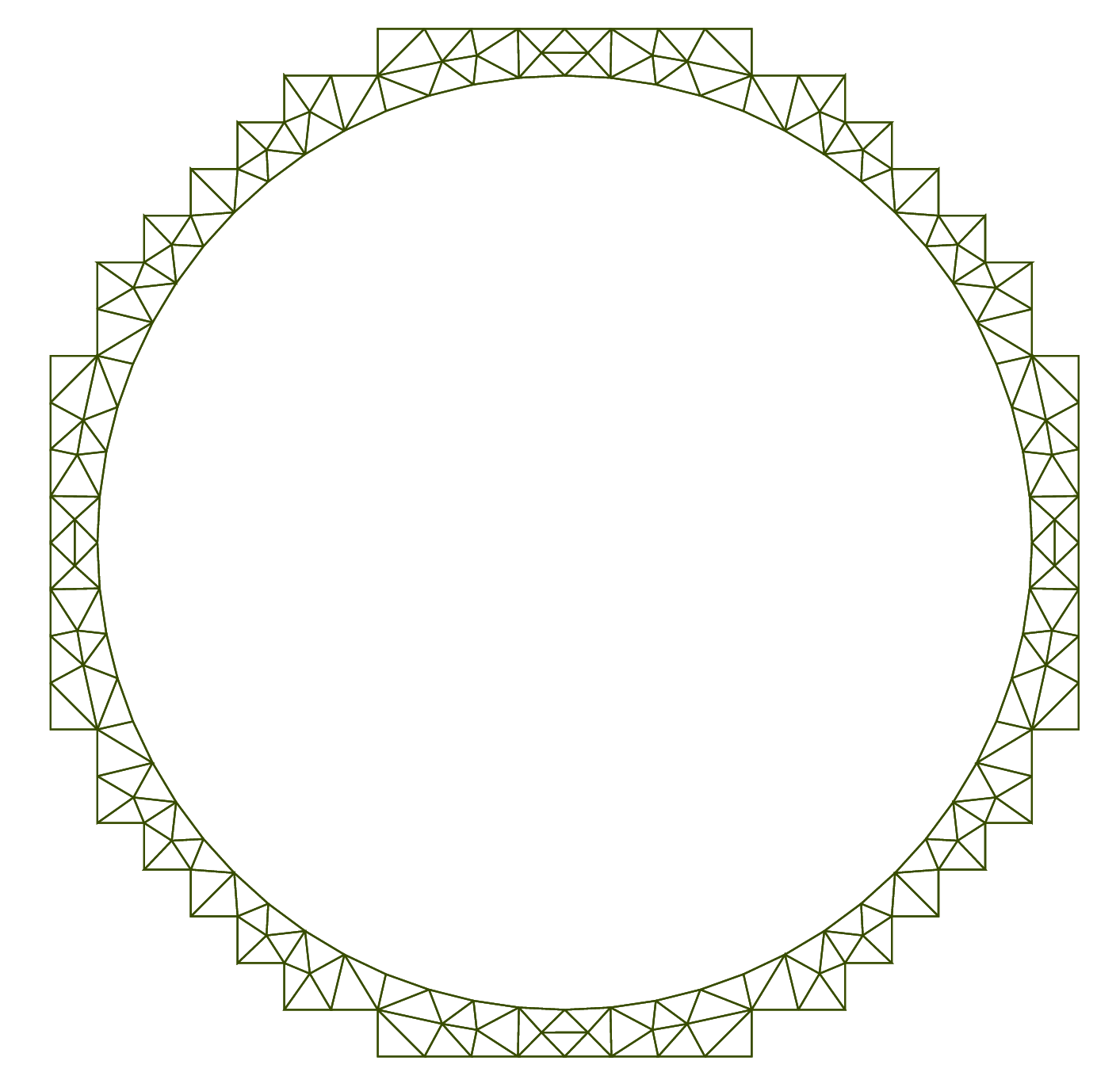}}
\caption{Scattering in a circular geometry: (a) Notation; (b) Minimal mesh for a circle of radius $R = 10h$ with the FEM mesh following a grid.}
\label{fig:TheMesh}
\end{figure}

The efficiency of the numerical method is illustrated in a two-dimensional problem with circle-shaped scatterer (see \fref{fig:TheMesh}(a)). It has analytical solution that is used as reference. The scatterer has a curved surface. This is important in order to demonstrate the improvement caused by the better geometry description of the FEM layer (coupled model) with respect to a staircase approximation based on the regular grid (use of only BAE \cite{poblet-PVS:2015}).

The force, which represents the imposed normal derivative of the variable $u$ at the contour, is chosen in order to generate a scattered wave described by means of only one cylindrical harmonic. The nodal values of the force vector are
\begin{equation}
 f_{i} = \cos\left( N \varphi_{i} \right), \qquad i \in \gamma_\rin
\end{equation}
The angle $\varphi$ and the radius $R$ of the circle are shown in the sketch of \fref{fig:TheMesh}(a). $N$ is related with the spatial frequency of the imposed force,  $N$ waves exist over the circle. The expression of the scattered field on the circle surface is
\begin{equation}\label{eq:AnalyticalSolution}
 u(R,\varphi) = \frac{2 H_{N}^{(1)}\left( \K R\right)}{H_{N-1}^{(1)}\left( \K R\right) - H_{N+1}^{(1)}\left( \K R\right)} \cos\left( N \varphi \right)
\end{equation}
where $H_{N}^{(1)}$ is the Hankel function of the first kind and order $N$ and $\K$ is the wavenumber of the problem.

Different error types play an important role in the numerical solution of this problem: \textit{i)}interpolation and dispersion error of the scattered field; \textit{ii)}error in the description of the oscillatory force imposed on the scatterer surface; and \textit{iii)}geometry error in the approximation of the scatterer shape. Error types \textit{i)} and \textit{ii)} are related with the number of nodes per wave length of the scattered field or the imposed force, respectively. Error type \textit{iii)} is related with the curvature of the scatterer. Each error type can be the dominant error source depending on the frequency range and the geometrical or material parameters of the model.

The mesh in \fref{fig:TheMesh}(b) is designed in order to have a transition zone between the circle (boundary $\Gamma_\rin$) and a closed grid shape. It is forced to be thin in order to use the minimum number of finite elements. This mesh has nodes $\omega_\rin$ and elements $\omega'_\rin$. The nodes over the internal boundary $\gamma_\rin$ are placed exactly on the circle (equally distributed). The force vector is null for nodes not belonging to $\gamma_\rin$. The nodes on the external boundary $\gamma_\rex$ are considered in the BAE part of the problem. The mesh is built with the GMSH software \cite{Geuzaine-Remacle:2009}.

The error is measured as
\begin{equation}
 e = \frac{\left|\left| \mathbf{u}^\rex_{\mathrm{num}} - \mathbf{u}^\rex_{\mathrm{exact}}\right|\right|}{\left|\left| \mathbf{u}^\rex_{\mathrm{exact}}\right|\right|} \simeq \sqrt{\frac{\sum_{i \in \gamma_\rex}^{n} \left| u^\rex_{\mathrm{num},i} - u^\rex_{\mathrm{exact},i}\right|}{\sum_{i \in \gamma_\rex}^{n} \left|u^\rex_{\mathrm{exact},i}\right|}}
\end{equation}
where `num' is the numerical solution and `exact' the solution obtained with (\ref{eq:AnalyticalSolution}).

In all the examples the grid spacing is $h = 1$. If nothing else is specified, the mean finite element size is also $h = 1$ and the layer of finite elements that surrounds the circle has an approximate external radius of $R_{ext} = R + h$.

\begin{figure}[ht]
\subfigure[]{\includegraphics[width=0.45\textwidth]{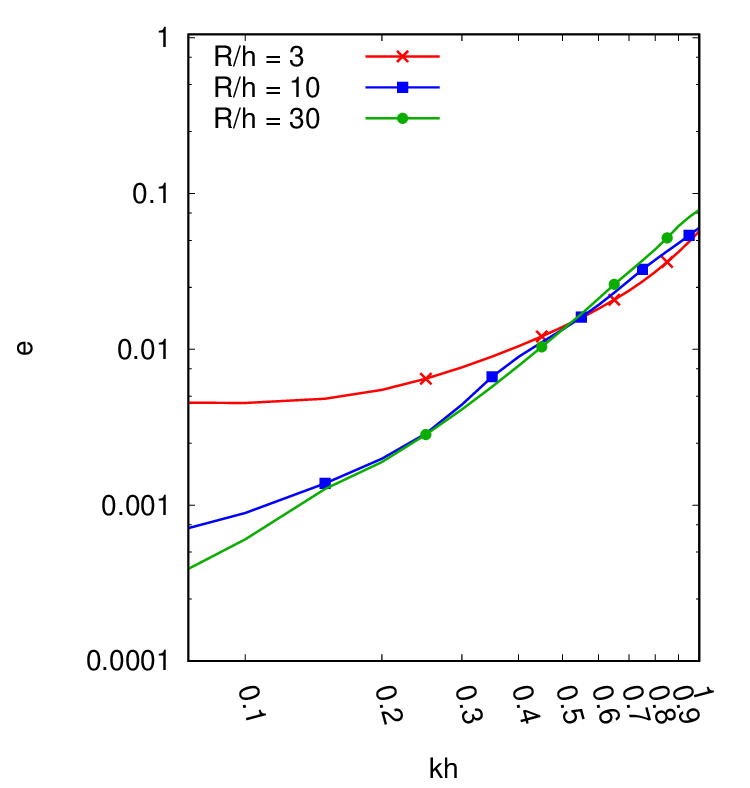}} \subfigure[]{\includegraphics[width=0.45\textwidth]{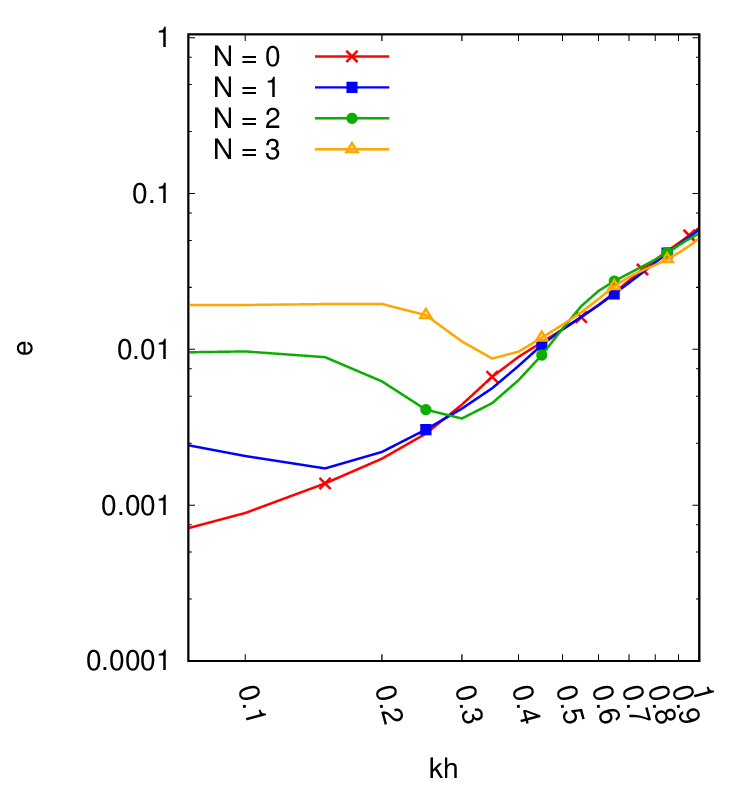}}
\caption{Relative errorfor different: (a) radius of the scatterer; (b) harmonics (shape of the imposed normal derivative).}
\label{fig:InfluenceOfRadiusAndHarmonic}
\end{figure}

\fref{fig:InfluenceOfRadiusAndHarmonic}(a) shows the error evolution with respect to the  dimensionless wavenumber  $\K h$ for several cylindrical scatterers of different size and the harmonic $N=0$. In all the cases the slope of the error curve is close to $2$. This is the expected result for the interpolation error of linear finite elements where $e = \theta h^{2}$, with $\theta$ a constant value \cite{Bouillard-Ihlenburg:99}.

It is observed that the numerical error has a different lower bound for each curve. This value is larger for smaller scatterers (with a more pronounced curvature compared to the element size) due to the geometry error of the linear finite element approximation of the circular shape. This error is invariant with respect to the wavenumber of the problem because it only depends on the relationship between the element size and the curvature of the circle. 

The geometry error is comparatively not important for large values of dimensionless wave number ($\K h \approx 0.3-1.0$) where the interpolation and dispersion error of the scattered field is dominant. On the contrary, geometry error becomes dominant at low frequencies when the scattered field is oscillating with a larger spatial wave length. As an example, consider the circle of radius $R = 3h$ where the exact curved piece of surface that contributes to each node is $ds \simeq \hfrac{2\pi R}{n} = 0.94247781$ ($n=20$ elements around the circle $\Gamma_\rin$). Its equivalent finite element length is $0.93860679$ which is slightly different.

For all this, it can be seen in  \fref{fig:InfluenceOfRadiusAndHarmonic}(a) how the theoretical convergence slope is lost for $\K h < 0.4$ in the circle of radius $R = 3h$ and for $\K h < 0.15$ in the circle of radius $R = 10h$. The circle of radius $R = 30h$ is not sensitive to the geometrical error in the studied frequency range.

The influence of the spatial wavenumber of the imposed force for a scatterer of radius $R = 10 h$ is shown in \fref{fig:InfluenceOfRadiusAndHarmonic}(b). There are $64$ nodes on the circle. The imposed force describes $N$ complete waves around the circle. Consequently, there are: $64$, $32$ and $16$ nodes per excitation wave length in the harmonics  $N=1,2,3$ respectively. This amount of nodes is related with the precision in the computation of the force vector. 

In the results of \fref{fig:InfluenceOfRadiusAndHarmonic}(b) two different zones can be clearly distinguished: large wavenumbers where the interpolation and dispersion error in $u^{\rex}$ is dominant and low frequencies where the error due to the force description is more important. Each curve has a limit wavenumber $\K$ for which the error in the solution becomes more or less constant and cannot be reduced with a decrease of $\K h$.  This limit value of the wavenumber $\K$ is related with the number of the harmonic $N$: $\K h \simeq 0.4$ for $N=3$, $\K h \simeq 0.3$ for $N=2$,  and $\K h \simeq 0.18$ for $N=1$. The curve corresponding to $N = 0$ is not affected by the error in the description of the force because it is constant all around the scatterer.

\begin{figure}[ht]
\subfigure[]{\includegraphics[width=0.45\textwidth]{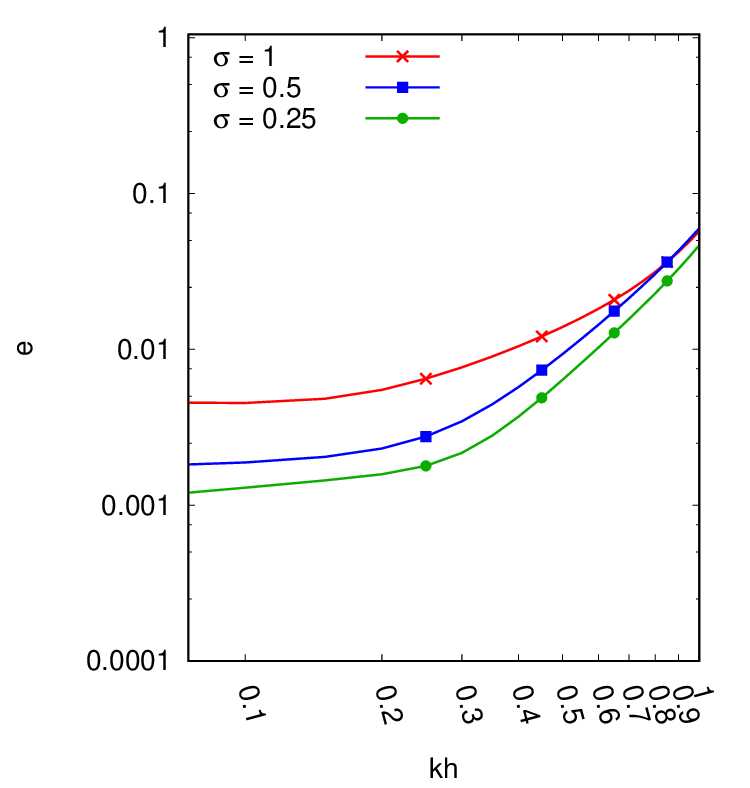}} \subfigure[]{\includegraphics[width=0.45\textwidth]{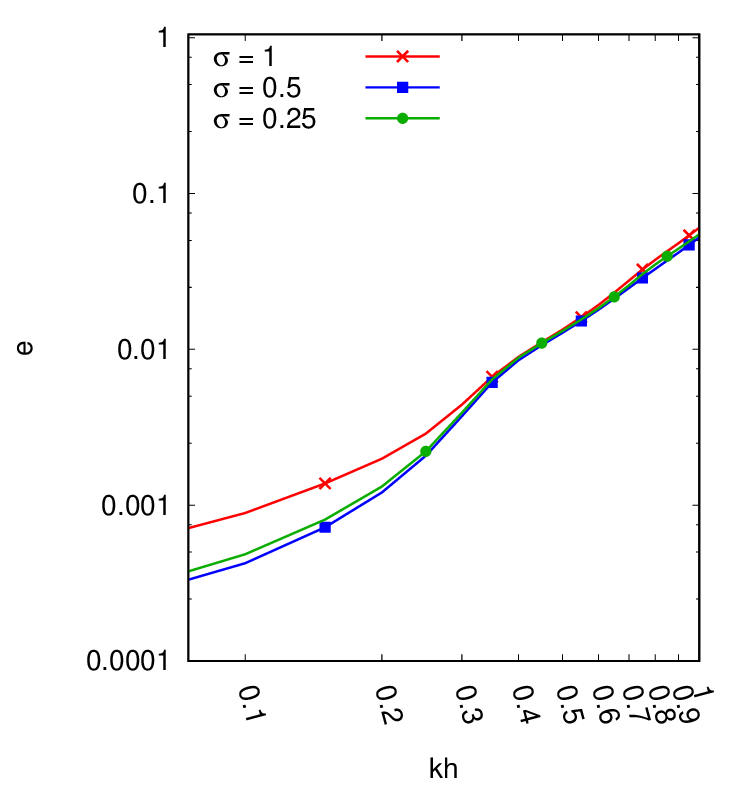}}
\caption{Influence of the finite element size around the scatterer for two different radius of the scatterer: (a) $R = 3h$; (b) $R = 10h$.}
\label{fig:InfluenceOfFemSize}
\end{figure}

\fref{fig:InfluenceOfFemSize} illustrates which is the effect of reducing the finite element size only on the circle (increase the number of nodes in $\gamma_\rin$) and not on the BAE contour (the number of nodes on $\gamma_\rex$ remains constant). The element size on $\Gamma_\rin$ is $\sigma h$, with $\sigma = 0.25, 0.5$ and $1$. The results are shown for two circles with radius $R = 3h$ and $R = 10h$. The improvement is more important for the case $R = 3h$ which is more sensitive to the geometry error at small wavenumbers. The reduction of the finite element size around the scatterer reduces the error in the whole frequency range. However, a lower bound (frequency invariant) is found for each $\sigma$ which shows again that it is due to approximation of the scatterer geometry and not due to the proper interpolation of the scattered field.

\begin{figure}[ht]
\subfigure[]{\includegraphics[width=0.45\textwidth]{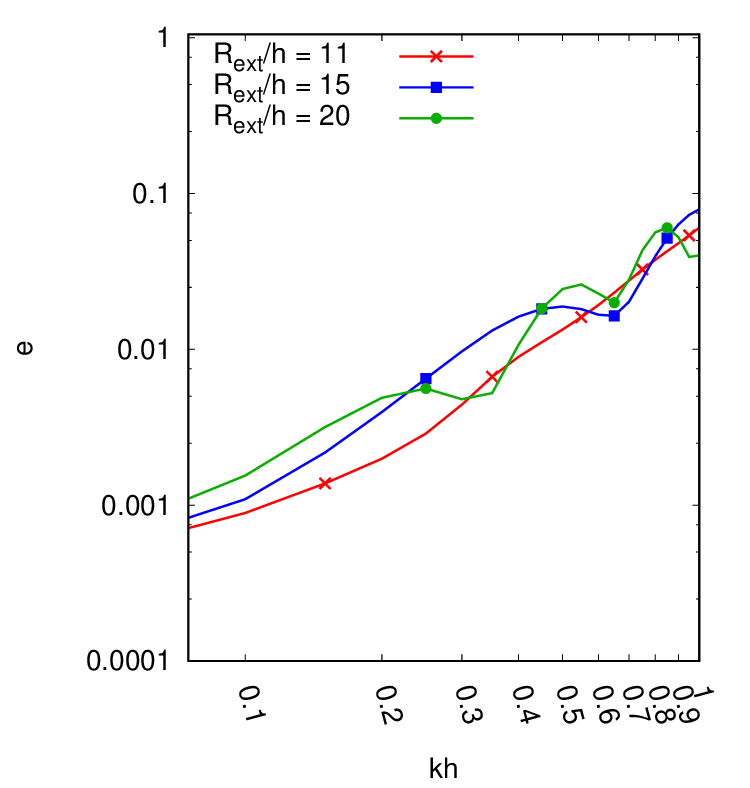}} \subfigure[]{\includegraphics[width=0.45\textwidth]{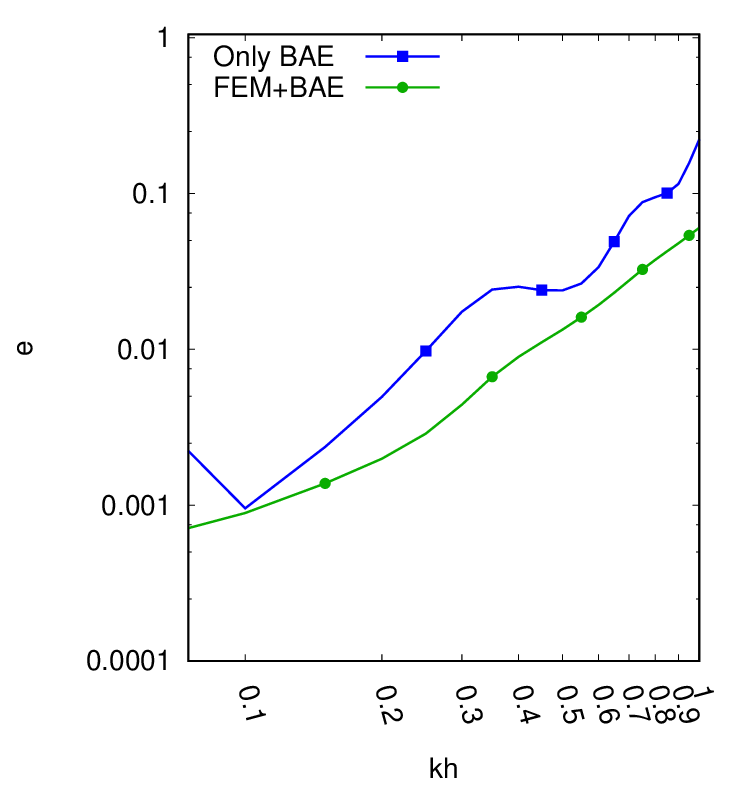}}
\caption{Influence of the domain truncation in the quality of the solution for a scatterer with raius $R = 10h$: (a) error for different FEM meshes truncated at $R_{ext}$; (b) comparison between a case where only BAE is used and the case when a small layer of finite elements is placed around the scatterer.}
\label{fig:InfluenceOfBoxSize}
\end{figure}

\fref{fig:InfluenceOfBoxSize} shows the effect of the finite element mesh truncation. First, the boundary $\Gamma^{\rex}$ is placed at several distances: $R+h$, $R+5h$ and $R+10h$ with a circular scatterer of radius $R = 10h$. It can be seen how the results are almost insensitive (or without clear meaningful trend) to the truncation distance. This is important because it allows the use of the thinnest finite element mesh around the scatterer, only conditioned by scatterer shape and meshing procedures. The use of a small mesh contributes to the reduction of computational costs. On the one hand, there are less unknowns. On the other hand, the range of required values of the discrete Green's function is smaller.

\fref{fig:InfluenceOfBoxSize}(b) shows a comparison between the case when the finite element layer is used (`FEM+BAE') and the case when it is not considered (`Only BAE'). In this second case the circular shape of the scatterer is approximated by means of a staircased geometry, defined by the closest grid (as it was done in \cite{poblet-PVS:2015}). 
One can observe the improvement caused by the description of the scatterer geometry by means of triangular finite elements comparatively to a grid approximation of the circle. The difference is larger for higher wavenumbers. But the slope or general trend is similar.

\section{Conclusions}
\label{sec:conclusions}

A numerical technique to deal with scattering problems has been presented. On the one hand, it can be understood as a complement to the CFIE--BAE method where a FEM layer is placed around the scatterer in order to better approximate its shape and reduce the geometry error. On the other hand, it can be understood as the use of BAE in order to exactly impose the radiation boundary condition in a FEM model.

It is shown how the resulting method keeps the properties of finite elements. Since linear triangles are considered for the FEM layer, order two convergence is observed. This behaviour is only truncated at very low values of dimensionless wavenumber $\K h$ by the geometrical error in the discretisation of the scatterer shape or the approximation of the force vector.

The coupling with FEM largely reduces the numerical error of BAE solutions and helps to overcome its main drawback in problems involving curve-shaped scatterers. This was caused by the staircase approximation of that shapes. That shapes are now approximated by means of standard finite elements without loosing any of the good properties of BAE for scattering problems: no need to compute boundary integrals (which are usually singular in other methods such as BEM), non-singularity of the problem even for the spurious eigenfrequencies of the scatterer and exact representation of the domain truncation.


\appendix
\section{Appendix A. Derivation of equations (\ref{eq0312a}) and (\ref{eq0308b})}
\label{sec:AppA}

First, derive (\ref{eq0312a}). Formally the proof can be written as follows.
Consider the expression
\[
\sum_{j,k \in \omega_\rex} u_j^\rex \alpha_{j,k}^\rex G_{k,m}
\]
On the one hand, due to (\ref{eq0302})
\begin{equation}
\sum_{j,k \in \omega_\rex} u_j^\rex \alpha_{j,k}^\rex G_{k,m} = \sum_{k \in \omega_\rex} h_k G_{k, m}.
\label{eqA01}
\end{equation}
On the other hand, due to (\ref{eq0308a}),
\begin{equation}
\sum_{j,k \in \omega_\rex} u_j^\rex \alpha_{j,k}^\rex G_{k,m} =
\sum_{j \in \omega_\rex} u_j^{\rex} (\delta_{j, m} + b_{j, m})
\label{eqA02}
\end{equation}
If $m \in \omega_\rex$, combining the expressions \ref{eqA01} and \ref{eqA02}, obtain
\begin{equation}
u_m^\rex = \sum_{j \in \gamma_\rex} (h_j^\rex G_{j,m} - u_j^\rex b_{j,m}).
\label{eqA03}
\end{equation}
After substitution  $u^\rex_j = u^\rin_j$ for $j \in \omega_\rex$ get (\ref{eq0312a}). Note that
(\ref{eqA03}) is valid only for the solution $u^\rex_j$ obeying the radiation condition.

However, this method cannot be applied directly, since the
summation is held over an infinite set of nodes $\omega_\rex$. In \cite{poblet-PVS:2014} one can find a refined procedure.
One should truncate the area $\Omega_\rex$, say, by a large square/cube, and apply (\ref{eqA01}), (\ref{eqA02}) to the truncated mesh. Then one should consider the limit of the size of the square/cube growing to infinity. The radiation condition obeyed by $u_j^\rex$ and $G_{j,m}$ guarantee that the integral over the outer boundary of the
sphere/cube vanishes.

Now apply matrix $\alpha^\rex_{m,n}$ to (\ref{eqA03}):
\begin{equation}
\sum_{m \in \omega_\rex} u_m^\rex \alpha^\rex_{m,n} =
\sum_{m \in \omega_\rex}\sum_{j \in \gamma_\rex} (h_j^\rex G_{j,m} - u_j^\rex b_{j,m}) \alpha^\rex_{m,n}.
\label{eqA04}
\end{equation}
Here the summation over $m$ causes no problem, since for each $n$ it is held only over the neighbors of~$n$,
where the coefficients $\alpha^\rex_{m,n}$ are non-zero.
Changing the order of summation in (\ref{eqA04}) and taking into account (\ref{eq0302}), get
\begin{equation}
h^\rex_n =
\sum_{j \in \gamma_\rex} \left( h_j^\rex \sum_{m \in \omega_\rex}G_{j,m} \alpha^\rex_{m,n} -
u_j^\rex \sum_{m \in \omega_\rex} b_{j,m} \alpha^\rex_{m,n} \right) .
\label{eqA05}
\end{equation}
Now multiply (\ref{eqA05}) by an arbitrary complex number $\nu$ with a non-zero imaginary part and add to
(\ref{eqA03}). The result is (\ref{eq0308b}).


\section{Appendix B. On invertibility of $\C$}
\label{sec:AppB}

The invertibility of $\C$ depends on details of realization of the finite element method, so here we can prove
a general but relatively weak theorem:

{\em If a homogeneous Dirichlet problem on $\Omega_\rex$ has no non-trivial solutions, then matrix $\C$
is invertible.}

A homogeneous Dirichlet problem on $\Omega_\rex$ is as follows: Find a function $w_j$ obeying equation
\begin{equation}
\sum_{j \in \omega_\rex}  \alpha^\rex_{m,j} w_{j} = 0,
\qquad
m \in (\omega_\rex \setminus \gamma_\rex),
\label{eqB01}
\end{equation}
boundary condition
\begin{equation}
w_j = 0 ,
\qquad
j \in \gamma_\rex,
\label{eqB02}
\end{equation}
and the radiation condition.

The uniqueness of solution of a homogeneous Dirichlet problem can be proven
in many particular cases.

The proof of the theorem is analogous to that of \cite{poblet-PVS:2014}.
Assume that all coefficients $\beta_{m,n}$ and $\alpha_{j,m}^\rex$ are real.
Let matrix $\C$ be not invertible.
This means that there exists a non-zero vector $v_j$, $j \in \gamma_\rex$ such that $v \C$ is a zero vector,
i.\ e.
\begin{equation}
\sum_{j \in \gamma_\rex} v_j G_{j,m} = \nu \sum_{j \in \gamma_\rex} v_j
\left(
\delta_{j,m} - G_{j,m} \alpha_{j,m}^\rex
\right)
,
\qquad
m \in \gamma_\rex.
\label{eqB03}
\end{equation}
Consider function $v_j$ on  $\gamma_\rex$. Introduce a ``single-layer potential''
on the uniform mesh $\omega$:
\begin{equation}
w_m = \sum_{j \in \gamma_\rex} G_{m,j} v_j,
\qquad
m \in \omega.
\label{eqB04}
\end{equation}
This function obeys equation (\ref{eqB01}) and the radiation condition by construction.
Note that
\begin{equation}
v_m = \sum_{j \in \bar \omega} G_{m,j} w_j.
\label{eqB05}
\end{equation}
Thus, (\ref{eqB03}) can be written in the form
\begin{equation}
w_m = \nu \sum_{j \in \bar \omega_\rin} \beta^\rin_{m,j} w_j,
\qquad
m \in \gamma_\rex,
\label{eqB06}
\end{equation}
where
\begin{equation}
\beta_{m,n}^\rin =
\left\{ \begin{array}{ll}
\beta_{m,n} & \mbox{if } m \mbox{ or } n \mbox{ belongs to }\bar \omega_\rin \setminus \gamma_\rex \\
\beta_{m,n} - \alpha_{m,n}^\rex & \mbox{otherwise}
\end{array} \right.
\label{eqB07}
\end{equation}
Note that
$\beta^\rin_{m,n} \ne 0$ only if $m,n \in \bar \omega_\rin$.
Note also that
\begin{equation}
 \sum_{j \in \bar \omega_in} \beta^\rin_{m,j} w_j =
 \sum_{j \in \bar \omega_in} \beta_{m,j} w_j = 0
\qquad
m \in (\bar \omega_\rin \setminus \gamma_\rex).
\label{eqB08}
\end{equation}

Consider a combination
\[
\sum_{m,n \in \bar \omega_\rin}
w^*_m \beta^\rin_{m,n} w_n
\]
where $\cdot^*$ denotes complex conjugation. Using (\ref{eqB06}) and (\ref{eqB08}) one can obtain
two representations for this combinations:
\begin{equation}
\sum_{m,n \in \bar \omega_\rin}
w^*_m \beta^\rin_{m,n} w_n =
\nu^{-1} \sum_{m \in \gamma_\rex}  w^*_m w_m =
(\nu^*)^{-1} \sum_{m \in \gamma_\rex} w^*_m w_m.
\label{eq09}
\end{equation}
Thus, we can conclude that
\begin{equation}
w_j = 0 ,
\qquad
j \in \gamma_\rex,
\label{eqB10}
\end{equation}
and $w_j$ is a solution of the homogeneous Dirichlet problem. It is non-trivial on $\omega_\rex$,
since equations (\ref{eqB05}) and
\begin{equation}
\sum_{j \in \bar \omega_\rin} \beta^\rin_{m,j} w_j =0 ,
\qquad
m \in \gamma_\rex,
\label{eqB11}
\end{equation}
(following from (\ref{eqB06})), are valid.

\section*{Acknowledgements}
The authors acknowledge the Euro-Russian Academic Network-Plus program (grant number 2012-2734/001-001-EMA2). J. Poblet-Puig from the LaC\`{a}N research group is grateful for the sponsorship/funding received from Generalitat de Catalunya (Grant number 2014-SGR-1471).
A.V.Shanin has been also supported by Russian Scientific school grant 7062.2016.2 and the Russian Foundation for Basic Research grant 14-02-00573.

\bibliographystyle{plain}
\bibliography{bibLegoBEM}

%
%
%
%

\end{document}